\documentclass[12pt,reqno]{amsart}

\usepackage[mathlines]{lineno}

\usepackage{times}
\usepackage{amsfonts}
\usepackage{amssymb}
\usepackage{latexsym}
\usepackage{xcolor}

\usepackage{comment}

\newtheorem{theorem}{Theorem}
\newtheorem{lemma}[theorem]{Lemma}

\newtheorem{remark}{Remark}
\newtheorem{corollary}[theorem]{Corollary}
\newtheorem{example}[theorem]{Example}

\def\nin{\relax\hbox{$/\kern-.7em{\rm \in\,}$}}

\begin{document}


\centerline{{\Large \bf Functions of  continuous Ces\'aro operators}}

\

\centerline{Adolf Mirotin }

\

\centerline{amirotin@yandex.ru}

\

ABSTRACT.  We describe  holomorphic functions and fractional powers  of Ces\'{a}ro operators in $L^2(\mathbb{R})$,  $L^2(\mathbb{R}_+)$, and $L^2[0,1]$. Logarithms of  Ces\'{a}ro operators are introduced as well and their spectral properties are studied. Several examples are considered.

\

\

\section{Introduction}

The study of discrete  Ces\'aro operator in $\ell^2$ and continuous Ces\'aro operators in  $L^2(\mathbb{R}_+)$  and $L^2[0,1]$ was initiated by  A. Brown, P.R. Halmos, and A.L. Shields \cite{BHSh}. On the other hand, during last two decades the theory of Hausdorff operators, the wast generalization of  Ces\'aro operators was developed (see the survey article \cite{LM:survey}).

It follows from the work of Wallen and Shields \cite{ShieldsWallen} that each discrete Hausdorff operator on $\ell^2$ is a bounded
holomorphic function of the discrete Ces\'aro operator on $\ell^2$ (see also \cite{KT}, \cite[Remark 2]{Leibowitz}).

  The continuous Ces\'aro operator originated from the classical summation methods for integrals   \cite[p. 110--111]{Hardy}, see also \cite{Boyd1, Boyd, Boyd2, Leibowitz2, Leibowitz3}. In the last  decade  continuous Ces\'aro operators continue to attract attention of researchers (see, e.g., \cite{ABR}, \cite{ABR1}, \cite{DLS} and the bibliography therein).

 In this work, we are aimed to describe  holomorphic functions of continuous Ces\'aro operators in $L^2(\mathbb{R})$,  $L^2(\mathbb{R}_+)$,  and $L^2[0,1]$. It turns out that the cases of the spaces $L^2(\mathbb{R}_+)$ and $L^2[0,1]$ are easily reduced to the case of $L^2(\mathbb{R})$.
 We also compute fractional powers of this operator  in  $L^2(\mathbb{R})$ and  $L^2(\mathbb{R}_+)$. In both cases, the result is the well-known H$\ddot o$lder operator.
 Logarithms of  Ces\'{a}ro operators are introduced and their spectral properties and converse are studied, as well. 

We show in particular that holomorphic functions of classical Ces\'aro operators, and also resolvents and converse of their logarithms are a sort of  Hausdorff operators.
Thus, more general  Hausdorff operators naturally arise as functions of continuous Ces\'aro operators. 
In this context, the paper also considers examples of several classical integral operators of harmonic analysis.

 Our considerations are based on the results of \cite{faa},  \cite{Forum}, and  \cite{LiflyandMirotin1}.

\section{Preliminaries}

A continuous Ces\'aro operator has the form

\begin{equation}\label{C}
f\mapsto \frac1x\int_{0}^xf(t)dt=\int_{0}^1f(ux)du.
\end{equation}

This is a very special case of a one-dimensional Hausdorff operator
\begin{equation}\label{HKa}
(\mathcal{H}_{K,a}f)(x):=\int_{\mathbb{R}}K(u)f(a(u)x)du,
\end{equation}
where   $K$ and $a$ are given measurable functions on $\mathbb{R}$.
If $a(u)=u$ on the support of $K$, we shall wright $\mathcal{H}_{K}$ instead of $\mathcal{H}_{K,a}$.

If the kernel $K$ satisfies the condition $K(u)|a(u)|^{-\frac12}\in L^1(\mathbb{R})$, then the operator $\mathcal{H}_{K,a}$ is bounded (by the Cauchy-Schwartz inequality) and normal in $L^2(\mathbb{R})$ and $L^2(\mathbb{R})_+$ (see, e.~g., \cite{Forum}). 

Note that the kernel $K$ is uniquely determined by the operator $\mathcal{H}_{K}$. Indeed, let
$$
\int_{\mathbb{R}}K(u)f(ux)du=\int_{\mathbb{R}}K_1(u)f(ux)du
$$
for all $f\in L^2(\mathbb{R})$,  $x\in\mathbb{R}$. Then, putting $f_{\pm}(t)=e^{-p|t|}\chi_{\mathbb{R}_{\pm}}(t)$ (${\rm Re}(p)>0$), where $\chi$ stands for the indicator function of the set indicated as a subscript, and $x=1$ in this equality we get $K=K_1$ a~.e.
by the uniqueness theorem for the Laplace transform.

The similar uniqueness result  is true for the restriction of  $\mathcal{H}_K$ to $L^2(\mathbb{R}_+)$.

We shall use the following facts from \cite{faa} (see also \cite{Forum}). Let $K(u)|a(u)|^{-\frac12}\in L^1(\mathbb{R})$. The matrix valued function
\begin{align}\label{phi3}
\Phi_{\mathcal{H}_K}(s)
&= \begin{pmatrix} \varphi_+(s)&\varphi_-(s)\\ \varphi_-(s)&\varphi_+(s) 
 \end{pmatrix} \ (s\in \mathbb{R})
\end{align}
where
\begin{align}\label{scalar}
\varphi_+(s)=\int_{(0,+\infty)}K(u)u^{-\frac12-\imath s}du,\\
 \varphi_-(s)=\int_{(-\infty,0)}K(u)|u|^{-\frac12-\imath s}du \nonumber
\end{align}
is called {\it  the matrix symbol} of a Hausdorff operator \eqref{HKa} (here and below $|u|^{-\imath s}=\exp(-\imath s\log |u|)$).

According to \cite{faa} we call the function
\begin{align}\label{symbol}
   \varphi(s):=\int_{\mathbb{R}}K(u)|u|^{-1/2- \imath s}du\ (s\in \mathbb{R})
\end{align}
\textit{ the scalar symbol} of a Hausdorff operator $\mathcal{H}_K$ in $L^2(\mathbb{R})$.
 
Next, in \cite{LiflyandMirotin1} the following set of Hausdorff operators was introduced 
$$
\mathcal{A}_a:=\left\{\mathcal{H}_{K,a}:  \frac{K(u)}{|a(u)|^{\frac{1}{2}}}\in L^1(\mathbb{R})\right\},
$$
for a fixed measurable  function $a$ on $\mathbb{R}$, $a(u)\ne 0$ for a.~e. $u\in \mathbb{R}$. It is shown that this set 
is a non-closed  commutative subalgebra of the algebra $\mathcal{L}(L^2(\mathbb{R}))$ of bounded operators on $L^2(\mathbb{R})$  without unit (for multidimensional analogs see \cite{LiflyandMirotin2}).
We shall write $\mathcal{A}_u$ instead of $\mathcal{A}_a$ if $a(u)=u.$

Moreover, in   \cite{LiflyandMirotin1} a holomorphic calculus for one-dimensional Hausdorff operators was considered and the next result was proved.
In the sequel,  $\sigma(T)$ denotes the spectrum of an operator $T$ in $L^2(\mathbb{R})$.

\begin{theorem}\label{th:3.4}(\cite[Theorem 3.1]{LiflyandMirotin1}) Let $\mathcal{H}_{K,a}\in \mathcal{A}_a$. If a function $F$ is holomorphic in the neighborhood  of the
 set  $\sigma(\mathcal{H}_{K,a})\cup\{0\}$ and $F(0)=0$, then $F(\mathcal{H}_{K,a})\in \mathcal{A}_a$.
\end{theorem}

\section{Functions of the Ces\'aro operator in $L^2(\mathbb{R})$}

In this section $\mathcal{C}$ stands for the Ces\'aro operator in $L^2(\mathbb{R})$. This operator is  normal,
   and its spectrum   equals to the circle $\mathbb{T}+1$ by \cite[Theorem 2]{faa} (here and below  $\mathbb{T}$ stands for the unit circle in the complex plane).

\subsection{Holomorphic functions of the Ces\'aro operator in $L^2(\mathbb{R})$}

Recall that the Mellin transform of a  function $f$ on $(0,+\infty)$ has the form
$$
(\mathcal{M}f)(z)=\int_{(0,+\infty)}f(u)u^{z-1}du.
$$

The next theorem  gives a full description for holomorphic functions of the Ces\'aro operator in $L^2(\mathbb{R})$.

\begin{theorem}\label{th:1} 1) Let $F$ be a holomorphic function in a neighborhood  of the spectrum $\mathbb{T}+1$ of the  Ces\'aro operator $\mathcal{C}$ in $L^2(\mathbb{R})$ and $F(0)=0$. Then $F(\mathcal{C})=\mathcal{H}_K$ for a unique  function $K$ and the following conditions hold:

(a) $K(u)=0$ for a.e. $u<0$,  

(b) $K(u)u^{-\frac12}\in L^1(\mathbb{R}_+)$, and 

(c) $(\mathcal{M}K)(z)=F(z^{-1})$ for all $z$ from the ``critical line'' $\mathrm{Re}(z)=\frac12$. 

2) If a function $K$ meets the conditions (a), (b), and (c), where $F$ is  holomorphic  in some neighborhood  of  $\mathbb{T}+1$, then  $F(0)=0$ and $\mathcal{H}_K=F(\mathcal{C})$. 
\end{theorem}

\begin{proof}
  1) Let $F$ be a holomorphic function in a neighborhood $N$ of  $\mathbb{T}+1$, $F(0)=0$. Then $F(\mathcal{C})\in \mathcal{A}_u$ by Theorem \ref{th:3.4}. In other wards, 
$F(\mathcal{C})=\mathcal{H}_K$, where $K$ meets the condition   $K(u)|u|^{-\frac12}\in L^1(\mathbb{R})$. The uniqueness of  $K$ was proved above. By \cite{faa} the matrix symbol of  $\mathcal{C}$ is

\begin{align}\label{matrixC}
\Phi_{\mathcal{C}}(s)=
\begin{pmatrix} (\frac12-\imath  s)^{-1}&0\\ 0& (\frac12-\imath s)^{-1}
\end{pmatrix},
\quad s\in \mathbb{R}, 
\end{align}
 since $\varphi_+(s)=\varphi_+^0(s):= (\frac12-\imath s)^{-1}$ and  $\varphi_-(s)=\varphi_-^0(s):= 0$ for the  Ces\'aro operator $\mathcal{C}$.

 Further, it is known  (see \cite{LiflyandMirotin1}), that
\begin{align*}
\Phi_{F(\mathcal{C})}
&= \begin{pmatrix} F_1(\varphi_-^0,\varphi_+^0)&F_2(\varphi_-^0,\varphi_+^0)\\ F_2(\varphi_-^0,\varphi_+^0)&F_1(\varphi_-^0,\varphi_+^0)
 \end{pmatrix}
\end{align*}
where
$$
F_1(z_1,z_2):=\frac{1}{2\pi i}\int_{\Gamma} F(\lambda)\frac{\lambda-z_2}{(\lambda-z_2)^2-z_1^2}d\lambda,
$$
and
$$
F_2(z_1,z_2):=\frac{-z_1}{2\pi i}\int_{\Gamma} \frac{F(\lambda)}{(\lambda-z_2)^2-z_1^2}d\lambda.
$$
Here $\Gamma$ denotes the boundary of any open annulus $U$ that contains  $\mathbb{T}+1$ such that $N$ contains its closure. Therefore,
\begin{align*}
\Phi_{F(\mathcal{C})}(s)
&= \begin{pmatrix} F_1(0,\varphi_+^0(s))&F_2(0,\varphi_+^0(s))\\ F_2(0,\varphi_+^0(s))&F_1(0,\varphi_+^0(s))
 \end{pmatrix}\\
&=\begin{pmatrix} F_1(0, (\frac12-\imath s)^{-1}&F_2(0, (\frac12-\imath s)^{-1}\\ F_2(0, (\frac12-\imath s)^{-1}&F_1(0, (\frac12-\imath s)^{-1})
 \end{pmatrix}.
 \end{align*}
Note that $F_20,z_2)=0$ and by the Cauchy's formula
$$
F_1(0,z_2)=\frac{1}{2\pi i}\int_{\Gamma} F(\lambda)\frac{\lambda-z_2}{(\lambda-z_2)^2}d\lambda= F(z_2),
$$
as $z_2\in U$.
This yields that for all $s\in \mathbb{R}$
\begin{align}\label{phi2}
\Phi_{F(\mathcal{C})}(s)
&= \begin{pmatrix} F((\frac12-\imath s)^{-1})&0\\ 0&F((\frac12-\imath s)^{-1})
 \end{pmatrix},
 \end{align}
because $(\frac12-\imath s)^{-1}\in \mathbb{T}+1$ as $s\in \mathbb{R}$.

On the one hand, formulas \eqref{phi3} and \eqref{scalar} hold.
Since $F(\mathcal{C})=\mathcal{H}_K$, we deduce firstly from \eqref{phi2} and  \eqref{phi3} that 
$$\int_{(-\infty,0)}K(u)|u|^{-\frac12-\imath s}du=0$$
  for all  $s\in \mathbb{R}$.
And secondly, 
$$
\int_{(0,+\infty)}K(u)u^{-\frac12-\imath s}du=F\left(\left(\frac12-\imath s\right)^{-1}\right)
$$
 for all  $s\in \mathbb{R}$. This implies (c).

2) Now let (a), (b), and (c) hold true for $K$ and $F$, where $F$ is a holomorphic function  in some neighborhood  of  $\mathbb{T}+1$. Since the Mellin transform
$$
(\mathcal{M}K)\left(\frac12+\imath t\right)=\int_0^\infty \left(K(u)u^{-\frac12}\right)u^{\imath t}du
$$
tends to zero as $t\to\infty$ by the Riemann-Lebesgue Lemma, we deduce that $F(0)=0$.

Next, it was shown in the first part of the proof that (a), (b), and (c) along with $F(0)=0$ imply that $\Phi_{F(\mathcal{C})}$ is given by the formula \eqref{phi2}.
On the other hand, since (a), we have by \eqref{phi3} that
\begin{align}\label{phi4}
\Phi_{\mathcal{H}_K}(s)
&= \begin{pmatrix} \varphi_+(s)&0\\ 0&\varphi_+(s)
 \end{pmatrix}
\end{align}
where
\begin{align*}
\varphi_+(s)&=\int_{(0,+\infty)}K(u)u^{-\frac12-\imath s}du\\
&=(\mathcal{M}K)\left(\frac12-\imath s\right)=F\left(\left(\frac12-\imath s\right)^{-1}\right)
\end{align*}
for all $s\in \mathbb{R}$. It follows that $\Phi_{\mathcal{H}_K}=\Phi_{F(\mathcal{C})}$ and thus $\mathcal{H}_K=F(\mathcal{C})$ (see, e.g., \cite[Lemma 2.1]{LiflyandMirotin1}).
\end{proof}

\begin{corollary}
Under the conditions of Theorem \ref{th:1} one has for the spectrum  $\sigma(\mathcal{H}_K)=F(\mathbb{T}+1)$.
\end{corollary}

\begin{remark}
Let the conditions (a), (b), and (c) hold. If $K$ is of bounded variation in a neighborhood of a point $x$, then by the inverse formula for Mellin transform
$$
\frac12(K(x+)+K(x-))=\frac{1}{2\pi\imath}v.p.\int_{\frac12-\imath\infty}^{\frac12+\imath\infty}F(s^{-1})x^{-s}ds.
$$

 \end{remark}

\begin{example}\label{ex:1}
Let ${\rm Re}(\alpha)<\frac12$. Consider the following operator in $L^2(\mathbb{R})$:
$$
(\mathcal{P}_\alpha f)(x)=|x|^{\alpha-1}{\rm sign}(x)\int_0^x|v|^{-\alpha} f(v)dv
$$
(Hardy operator). Then $\mathcal{P}_\alpha =\mathcal{H}_K$ where $K(u)= |u|^{-\alpha}\chi_{(0,1)}(u)$. All the conditions (a), (b), and (c) hold, because
$$
(\mathcal{M}K)(z)=\int_0^1u^{-\alpha+z-1}du=\frac{1}{z-\alpha}.
$$
In this case, 
\begin{equation}\label{Falpha}
F(z)=F_\alpha(z):=\frac{z}{1- \alpha z}
\end{equation}
  is  holomorphic in a neighborhood of $\mathbb{T}+1$, and so 
  $$
 \mathcal{P}_\alpha =F_\alpha(\mathcal{C})= \mathcal{C}\left(I-\alpha \mathcal{C}\right)^{-1}=\alpha^{-1}((I-\alpha\mathcal{C})^{-1}-I)
  $$
if  $\alpha\ne 0$.  Putting here $\lambda=\alpha^{-1}$   we get for the resolvent of the Ces\'{a}ro operator in $ L^2(\mathbb{R})$
\begin{equation}\label{R}
R(\lambda,\mathcal{C}):=(\lambda I-\mathcal{C})^{-1}=\lambda^{-2}P_{\lambda^{-1}}+\lambda^{-1} I
\end{equation}
for all $\lambda$ exterior to the circle  $\mathbb{T}+1$.

   It follows, that the spectrum  $\sigma(\mathcal{P}_\alpha)$ equals to the circle $F_\alpha(\mathbb{T}+1)$.
\end{example}

\begin{example}\label{ex:2}
Let ${\rm Re}(\alpha)<\frac12$. Consider the following operator in $L^2(\mathbb{R})$:
$$
(\mathcal{Q}_\alpha f)(x)=
\begin{cases}
x^{-\alpha}\int_x^\infty v^{\alpha-1}f(v)dv, & when \ x>0,\\
|x|^{-\alpha}\int_{-\infty}^x |v|^{\alpha-1}f(v)dv, & when \ x<0\\
\end{cases}
$$
($\mathcal{Q}_0$ is the so called Copson operator).
  Then $\mathcal{Q}_\alpha\in \mathcal{A}_u$ with
the kernel $K(u)=\chi_{(1,\infty)}(u)|u|^{\alpha-1}$.
The conditions (a) and (b) are obviously satisfied. The condition (c) holds too, because
$$
(\mathcal{M}K)(z)=\int_1^\infty u^{\alpha+z-2}du=\frac{1}{1-z-\alpha},
$$
if (and only if) ${\rm Re}(z+\alpha)-1< 0$.
Thus,  
\begin{equation}\label{F}
F(z)=\frac{ z}{(1-\alpha)z-1}
\end{equation}
  is  holomorphic in a neighborhood of $\mathbb{T}+1$, and so for $\alpha\ne 1$
  $$
\mathcal{Q}_\alpha=F(\mathcal{C})=\mathcal{C}((I-\alpha)\mathcal{C}-I)^{-1}=-(1-\alpha)^{-1}((I-(1-\alpha)\mathcal{C})^{-1}-I).
  $$
Putting here $\lambda=(1-\alpha)^{-1}$   we get for the resolvent of the Ces\'{a}ro operator in $L^2(\mathbb{R})$
$$
R(\lambda,\mathcal{C})=(\lambda I-\mathcal{C})^{-1}=-\lambda^{-2} \mathcal{Q}_{1-\lambda^{-1}}+ \lambda^{-1} I
$$
 for all $\lambda$ interior to the circle  $\mathbb{T}+1$.

Thus, we conclude that
$$
R(\lambda,\mathcal{C})=
\begin{cases}
\lambda^{-2}\mathcal{P}_{\lambda^{-1}}+\lambda^{-1} I, when \  |\lambda-1|>1\\
-\lambda^{-2} \mathcal{Q}_{1-\lambda^{-1}}+ \lambda^{-1} I,  when \  |\lambda-1|<1.
\end{cases}
$$
This is consistent with the result of D.~ Boyd \cite{Boyd2} for the resolvent of the Ces\'{a}ro operator in $L^2(\mathbb{R}_+)$.

  We conclude also that the spectrum  $\sigma(\mathcal{Q}_\alpha)$ equals to the circle $F(\mathbb{T}+1)$.

\end{example}

\begin{example}\label{ex:3}
Consider the following Hausdorff operator in $L^2(\mathbb{R})$:
$$
(Rf)(x)=\int_0^1\left\{\frac1u\right\}f(ux)du
$$
where $\{t\}:=t-\lfloor t \rfloor$ denotes the fractional part of a real number $t>0$. The kernel $K(u)=\left\{\frac1u\right\}\chi_{(0,1)}(u)$ of this operator satisfies  conditions (a) and (b). On the other hand,
it follows from the theory of the Riemann zeta function (see, e.g. \cite{T}) that
$$
\zeta(s)=\frac{s}{s-1}-s\int_1^\infty\{t\}t^{-1-s}dt,  \ {\rm Re}(s)>0, s\ne 1.
$$

Therefore,
$$
(\mathcal{M}K)(s)=\int_0^1\left\{\frac1u\right\}u^{s-1}du= \frac{1}{s-1}- \frac{1}{s}\zeta(s), \ {\rm Re}(s)>0, s\ne 1.
$$
Here the function
$$
F(z)=\frac{z}{1-z}-z\zeta\left(\frac{1}{z}\right)
$$
 is not holomorphic in any neighborhood of $\mathbb{T}+1$, since 
the zeta function has an essential singularity at  infinity   \cite{SA}. Thus,  the operator $R$     is not  a  holomorphic function of  $\mathcal{C}$.
\end{example}

\subsection{Fractional powers and the logarithm of the Ces\'aro operator in $L^2(\mathbb{R})$}

 In \cite{DLS} positive fractional powers  of the Ces\'aro operator in $L^\infty(\mathbb{R}_+)$  were defined and the   limits at infinity of the values of  fractional powers were considered.
In this section, we compute complex fractional powers  of the Ces\'aro operator in $L^2(\mathbb{R})$.  We consider the  branch  of the function $z^\alpha$ with $\mathrm{Re}(\alpha)>0$ in the right half-plane $\{\mathrm{Re}(z)>0\}$ given by the  condition $1^\alpha=1$ and put $0^\alpha:=0$. This function is  continuous (but not holomorphic)  on the spectrum $\mathbb{T}+1$  of  $\mathcal{C}$. Fractional powers  of the normal operator in $L^2(\mathbb{R})$ will be considered in the sense of the functional calculus for normal operators (see, e.g., \cite{DSch2}, or \cite{Conway}). 

 Recall that {\it the H$\ddot o$lder means  of order $\alpha$} is  defined as
$$
(H_\alpha f)(x)=\frac{1}{\Gamma(\alpha)}\int_0^1\left(\log \frac1u\right)^{\alpha-1}f(ux)du.
$$ 
This is a Hausdorff operator of the form $\mathcal{H}_K$ which  is bounded in $L^2(\mathbb{R})$ if $\mathrm{Re}(\alpha)>0$. For real  $\alpha$ it is well known to workers in summability theory. The study of this operator was initiated by Rhoades \cite{Rhoades}, see also \cite{Leibowitz2}.

\begin{theorem}\label{th:2} Let $\mathrm{Re}(\alpha)>0$. Then
\begin{align*}
\mathcal{C}^\alpha=H_\alpha.
\end{align*}
\end{theorem}

\begin{proof} According to \cite[Theorem 2]{faa} (or  \cite[Theorem 1]{Forum}) operator
$\mathcal{C}$ in $L^2(\mathbb{R})$ with matrix symbol $\Phi_{\mathcal{C}}$ is unitary equivalent to the
operator $M_{\Phi_{\mathcal{C}}}$ of multiplication by the normal matrix-function $\Phi_{\mathcal{C}}$ in the space $L^2(\mathbb{R};\mathbb{C}^{2})$ of $\mathbb{C}^{2}$-valued
functions. This means that 
$$
\mathcal{C}=\mathcal{V}^{-1}M_{\Phi_{\mathcal{C}}}\mathcal{V}
$$
for some unitary operator $\mathcal{V}$ between $L^2(\mathbb{R})$ and $L^2(\mathbb{R};\mathbb{C}^{2})$. In particular, $\sigma(\mathcal{C})=\sigma(M_{\Phi_{\mathcal{C}}})=:\sigma$. Let $E_{M_{\Phi}}$ be the spectral measure of $M_{\Phi_{\mathcal{C}}}$. Since, according to the spectral theorem,
$$
M_{\Phi_{\mathcal{C}}}=\int_\sigma \lambda dE_{M_{\Phi}}(\lambda),
$$
we have
$$
\mathcal{C}=\int_\sigma \lambda d(\mathcal{V}^{-1}E_{M_{\Phi}}(\lambda)\mathcal{V}).
$$
It follows, that  $\mathcal{V}^{-1}E_{M_{\Phi}}\mathcal{V}$ is the spectral measure for $\mathcal{C}$  (see, e.g., \cite[Chapter X, \S 2, Corollary 7]{DSch2}), and therefore
\begin{align*}
\mathcal{C}^\alpha&=\int_\sigma \lambda^{\alpha} d(\mathcal{V}^{-1}E_{M_{\Phi}}(\lambda)\mathcal{V})\\
&=\mathcal{V}^{-1}\left(\int_\sigma \lambda^{\alpha} d(E_{M_{\Phi}}(\lambda)\right)\mathcal{V}\\ 
&=\mathcal{V}^{-1}M_{\Phi}^{\alpha}\mathcal{V}\\
&=\mathcal{V}^{-1}M_{\Phi^{\alpha}}\mathcal{V},
\end{align*}
where   in view of \eqref{matrixC}
$$
\Phi^{\alpha}\equiv\Phi_\mathcal{C}^{\alpha}=((\varphi_+^0) I_2)^{\alpha}=
(\varphi_+^0)^\alpha I_2
$$  
($I_2$ stands for the unit matrix of order two).  So, 
$\Phi_{\mathcal{C}^{\alpha}}=(\varphi_+^0)^{\alpha} I_2$. 

On the other hand, the operator $H_\alpha$ has the kernel
$$
K(u)=\frac{1}{\Gamma(\alpha)}\left(\log \frac1u\right)^{\alpha-1}\chi_{(0,1)}(u),
$$
and thus its matrix symbol $\Phi_{H_\alpha}$ has the form \eqref{phi4}, where (see, e.g., \cite[p.38, 4.37]{Ober})
\begin{align}\label{symbolH}
\varphi_+(s)&=\int_{(0,+\infty)}K(u)u^{-\frac12-\imath s}du\\  \nonumber
&=\frac{1}{\Gamma(\alpha)}\int_0^1\left(\log\frac1u\right)^{\alpha-1}u^{-1/2-\imath s}du\\ \nonumber
&=\left(\frac 12-\imath s\right)^{-\alpha} \nonumber
=(\varphi_+^0)^{\alpha}(s). \nonumber
\end{align}
It follows that $\Phi_{H_\alpha}=\Phi_{\mathcal{C}^{\alpha}}$, which completes the proof.
\end{proof}

\begin{corollary}\label{cor:5}Let  $\mathrm{Re}(\alpha)>0$. Then
$$
\sigma(H_\alpha)=\sigma(\mathcal{C}^{\alpha})=\{z^{\alpha}:z\in \mathbb{T}+1\},  
$$
and
$$
\|H_\alpha\|=\|\mathcal{C}^{\alpha}\|=\left(\frac{2\mathrm{Re}(\alpha)}{|\alpha|}\right)^{\mathrm{Re}(\alpha)}e^{\mathrm{Im}(\alpha)\arg(\alpha)}.
$$
In particular,   if $\alpha>0$, then $\sigma(H_\alpha)=\sigma(\mathcal{C}^{\alpha})$ is an ark
$$
\left\{((2\cos t)^\alpha\cos\alpha t, (2\cos t)^\alpha\sin\alpha t): t\in \left[-\frac\pi 2, \frac\pi 2\right]\right\},  
$$
and
$$\|H_\alpha\|=\|\mathcal{C}^{\alpha}\|=2^{\alpha}.
$$

\end{corollary}

\begin{proof} The first statement follows from the spectral mapping theorem.  Further, since $H_\alpha$ is normal, we have 
$$
\|H_\alpha\|=\sup\{|z^{\alpha}|:z\in \mathbb{T}+1\}=\sup\{|z|^{\mathrm{Re}(\alpha)}e^{-\mathrm{Im}(\alpha)\arg(z)}: z\in \mathbb{T}+1\}.
$$
Let $\alpha=a+\imath b, a>0$, $z=1+e^{\imath \theta},  \theta\in [-\pi,\pi]$. Since $\arg(1+e^{\imath \theta})=\theta/2,$ and
$|1+e^{\imath \theta}|=2\cos \theta/2$, 
 one has
$$
\|H_\alpha\|=\max_{ \theta\in [-\pi,\pi]}\left(2\cos\frac \theta 2\right)^{a}e^{-b\frac \theta 2}=2^a\max_{ x\in [-\pi/2,\pi/2]}\left(\cos x\right)^{a}e^{-bx}.
$$
The function $\left(\cos x\right)^{a}e^{-bx}$ attains its maximum on $[-\pi/2,\pi/2]$ at the point $x=-\tan^{-1}(b/a)=-\arg(\alpha)$ and the statement for $\|H_\alpha\|$ follows.

The formulas for the special case $\alpha>0$ follow from the fact that if $z\in \mathbb{T}+1$ and $\arg z=t$, then $|z|=2\cos t/2$.
\end{proof}


\begin{corollary}\label{cor:7} The operator $H_\alpha$ is a holomorphic function of $\mathcal{C}$ if and only if $\alpha\in\mathbb{Z}_+$.
\end{corollary}

\begin{corollary}\label{cor:8}  If we put $H_0:= I$, then  $H_\alpha$ (as a family of operators in $L^2(\mathbb{R})$)  is a semigroup which is  holomorphic in the right half-plane  $\{\alpha:\mathrm{Re}(\alpha)>0\}$.  The generator  of this semigroup is a normal operator whose  spectrum  lies in some half-plane  $\{z:\mathrm{Re}(z)\le \gamma\}$ with $\gamma<\infty$. The restriction of this semigroup to  the semiaxis  $\{\alpha\ge 0\}$ is a  positive semigroup.
\end{corollary}
We shall see below that one can take $\gamma=\log 2$ and this value is sharp.
\begin{proof} Indeed, it was shown in the previous theorem that $H_\alpha$ is unitary equivalent to the multiplication semigroup $M_{(\varphi_+^0)^{\alpha}}$ which is a  holomorphic semigroup in the right half-plane. The next statement follows from \cite[Theorem 13.37]{Rudin}. The last statement is obvious. 
\end{proof}
It is natural to denote   the generator of the semigroup $H_\alpha=\mathcal{C}^\alpha$ by $\log\mathcal{C}$ (cf. \cite[Proposition 10.1.2, Theorems 10.1.3,  10.1.4]{Martinez-Sanz}).

\begin{corollary} The semigroup $T(t)=2^{-t}H_t$ is strongly stable, i.~e. 
$$
\lim_{t\to +\infty}2^{-t}\|H_t f\|=0
$$
for all  $f\in L^2(\mathbb{R})$. 
\end{corollary}

\begin{proof}  Since $\forall t\ge 0  \|T(t)\|=1$, the semigroup $T$  is   Ces\'aro ergodic in $L^2(\mathbb{R})$  (see, e.g., \cite[Corollary 4.3.5]{Arendt}). 
This means that  for each  $f\in L^2(\mathbb{R})$ the limit $f_1=\lim_{t\to +\infty}2^{-t}H_t f$ exists.
Further,  by \cite[Proposition 4.3.1]{Arendt} for each  $f\in L^2(\mathbb{R})$  vector $f_1$ belongs to the kernel $\mathrm{Ker}(A)$, where $A$ is the generator of $T$.
Assume that  $f_1\ne 0$ for some $f\in L^2(\mathbb{R})$. Since $A=\log \mathcal{C}-(\log 2) I$, it follows that $\log 2$ is an eigenvalue of $\log \mathcal{C}$. 
On the other hand, the point spectrum 
\begin{equation}\label{sigmap}
\sigma_p(\log \mathcal{C})=\varnothing,
\end{equation} 
and we get a contradiction. To prove the  equality \eqref{sigmap}, note that according to \cite[Theorem 1]{faa}, the point spectrum of $H_t$ equals to the  point spectrum  of the operator $M'_{\varphi_+}$ of multiplication by $\varphi_+$ in $L^2(\mathbb{R})$, where $\varphi_+(s)=(1/2-\imath s)^{-t}$ is the scalar  symbol of $H_t$ by \eqref{symbolH}. But  $\sigma_p(M'_{\varphi_+})=\varnothing$ and so, $\sigma_p(H_t)=\varnothing$, too. It remains to apply the spectral mapping theorem for the point spectrum that states that  $\sigma_p(H_t)\setminus\{0\}=\exp(t\sigma_p(\log \mathcal{C}))$ (see, e.g., \cite[p. 277]{EN}), and the result follows.
\end{proof}

In the following
$$
\nu(y,-1)=\int_{-1}^\infty\frac{y^s}{\Gamma(s+1)}ds=\int_{0}^\infty\frac{y^{t-1}}{\Gamma(t)}dt
$$
(see, e.g., \cite{BEIII}).

Recall that for any linear operator $A$ its {\it spectral bound}
is defined by
$$
s(A) := \sup\{\mathrm{Re}(\lambda): \lambda\in \sigma(A)\}.
$$

The next theorem describes  spectral properties of the operator $\log\mathcal{C}$.

\begin{theorem}\label{th:resolvent} The following assertions hold true.

(i)  Let $\lambda>\log 2$. Then the resolvent of  $\log\mathcal{C}$ is a Hausdorff operator of the form 
\begin{align}\label{res}
R(\lambda,\log\mathcal{C})f(x)=e^{-\lambda}\int_0^1 \nu(e^{-\lambda}\log \frac1u,-1) f(ux)du,\ f\in L^2(\mathbb{R}).
\end{align}

(ii)  Let $\lambda>\log 2$. Then the spectrum of $R(\lambda,\log\mathcal{C})$ is
$$
\sigma(R(\lambda,\log\mathcal{C}))=\left\{\frac{1}{\lambda+\log(\frac12-\imath s)}:s\in \mathbb{R}\right\}\cup\{0\}.
$$

(iii) The spectrum of $\log\mathcal{C}$ equals to its continuous spectrum and is a curve
\begin{align}\label{spec(log C)}
\sigma(\log\mathcal{C})=\left\{\log(2\cos y)+\imath y: y\in \left(-\frac\pi 2, \frac\pi 2\right)\right\}.
\end{align}
In particular, the spectral bound  $s(\log\mathcal{C})=\log 2$. 

(iv) Each $\lambda\in \mathbb{R}$, $\lambda\ne \log 2$ is a regular point for 
$\log\mathcal{C}$ and  we have for the norm of the resolvent
\begin{align}\label{norm}
\|R(\lambda,\log\mathcal{C})\|=\frac{1}{\lambda-\log 2} \mbox{ for } \lambda> \log 2, \lim_{\lambda\to -\infty} \|R(\lambda,\log\mathcal{C})\|=\frac2\pi.
\end{align}
\end{theorem}

\begin{proof} (i)  For $\lambda>\log 2$, $f\in L^2(\mathbb{R})$, and $x\in \mathbb{R}$  one has  by the Fubini's theorem 
\begin{align}\label{est10}
R(\lambda,\log\mathcal{C})f(x)&=\int_0^\infty e^{-\lambda\alpha}\mathcal{C}^\alpha f(x)d\alpha\\ \nonumber
&=\int_0^\infty e^{-\lambda\alpha}H_\alpha f(x)d\alpha\\ \nonumber
&=\int_0^\infty e^{-\lambda\alpha} \int_0^1\frac{1}{\Gamma(\alpha)}(\log\frac1u)^{\alpha-1}f(ux)dud\alpha\\ \nonumber
&=\int_0^1\left(\int_0^\infty e^{-\lambda\alpha} \frac{1}{\Gamma(\alpha)}(\log\frac1u)^{\alpha-1}d\alpha\right)f(ux)du\\ \nonumber
&=\int_0^1\left(\int_{-1}^\infty e^{-\lambda(s+1)} \frac{1}{\Gamma(s+1)}(\log\frac1u)^{s}ds\right)f(ux)du\\ \nonumber
&=e^{-\lambda}\int_0^1 \nu(e^{-\lambda}\log \frac1u,-1)f(ux)du.
\end{align}
To justify the application of the Fubini's theorem in \eqref{est10} and to guarantee the convergence of the last integral, note that for each integer  $N$
$$
\nu(y,-1)=\exp(y)+O(|y|^{-1-N}), y\to +\infty
$$
(see, e.g., \cite[\S 18.3]{BEIII}), and consider the function
\begin{align*}
g(u):&=\int_{-1}^\infty \left|e^{-\lambda s} \frac{1}{\Gamma(s+1)}(\log\frac1u)^{s}\right|ds\\ 
&=\int_{-1}^\infty e^{- s\lambda} \frac{1}{\Gamma(s+1)}(\log\frac1u)^{s}ds\\ 
&=\nu(e^{-\lambda}\log\frac1u,-1)\\ 
&=\exp(e^{-\lambda} \log\frac1u)+O\left(\left|e^{-\lambda} \log\frac1u \right|^{-1-N}\right) (u\to 0+) \\  
&=\frac{1}{u^{e^{-\lambda}}}+O\left(\left(e^{-\lambda} \log\frac1u\right)^{-1-N}\right) (u\to 0+).  
\end{align*}
Since $2e^{-\lambda}<1$, the first summand in the last expression belongs to $L^2(0,1)$. If $N=-1$, the second summand  belongs to $L^2(0,1)$, too. Thus, $g(u)\in L^2(0,1)$. Since $f(ux)\in L^2(0,1)$ for all $x$,  the application of Fubini's theorem in \eqref{est10} is correct.

(ii) This follows from \cite[Theorem 1]{faa} and the fact that by (i) for  $\lambda>\log 2$ the scalar symbol of $R(\lambda,\log\mathcal{C})$ is
\begin{align}\label{varphi}
\varphi(s)=\frac{1}{\lambda+\log(\frac12-\imath s)}.
\end{align}

To prove this equality, note that in the case under consideration the kernel
$$
K(u)=e^{_\lambda}\nu(e^{-\lambda}\log \frac1u,-1)\chi_{(0,1)}(u).
$$
Then according to \eqref{symbol}
\begin{align*}
\varphi(s)&=e^{-\lambda}\int_0^1\nu(e^{-\lambda}\log \frac1u,-1)u^{-\frac12-\imath s}du\\ 
&=e^{-\lambda}\int_0^1\int_{0}^\infty\frac{\left(e^{-\lambda}\log \frac1u\right)^{t-1}}{\Gamma(t)}dtu^{-\frac12-\imath s}du\\
&=\int_{0}^\infty\left(\int_0^1u^{-\frac12-\imath s}\left(\log \frac1u\right)^{t-1}du\right)\frac{e^{-\lambda t}}{\Gamma(t)}dt\\
&=\frac{1}{\lambda+\log(\frac12-\imath s)}.
\end{align*}
Above we used the formula
$$
\int_0^1x^{\alpha-1}\log^\sigma\frac1x dx=\alpha^{-\sigma-1}\Gamma(\sigma+1)
$$
which is valid for $\mathrm{Re}(\alpha)>0, \mathrm{Re}(\sigma)>-1$ (see, e.g., \cite[p. 38, 4.37]{Ober}).

(iii) Let $\lambda>\log 2$. It is well known that 
$$
\sigma(R(\lambda,\log\mathcal{C}))\setminus\{0\}=\left\{\frac{1}{\lambda-z}:z\in \sigma(\log\mathcal{C})\right\}
$$
(see, e.g., \cite[Theorem 1.13]{EN}). Comparing this with (ii) we have
$$
\left\{\frac{1}{\lambda+\log(\frac12-\imath s)}:s\in \mathbb{R}\right\}=\left\{\frac{1}{\lambda-z}:z\in \sigma(\mathcal{C})\right\}.
$$
Therefore the spectrum of $\log\mathcal{C}$ is 
the curve
\begin{align*}
\sigma(\log\mathcal{C})=\{-\log(\frac12-\imath s):s\in \mathbb{R}\}
\end{align*}
with the parametric equations 
$$
x=-\log\sqrt{\frac14+s^2}\ \mbox{ and  } y=\tan^{-1}(2s), \  s\in \mathbb{R}.
$$
It remains to note that $y\in(-\frac\pi 2, \frac\pi 2)$ and
$$
2\cos y=2\sqrt{\cos^2(\tan^{-1}(2s))}=\frac{2}{\sqrt{1+\tan^2(\tan^{-1}(2s))}}=\frac{1}{\sqrt{\frac14+s^2}}.
$$
Thus, the spectrum of $\log\mathcal{C}$ is the graph of the function  $x=\log(2\cos y)$ ($y\in(-\frac\pi 2, \frac\pi 2)$), and \eqref{spec(log C)} follows.

Finally, the point spectrum of $\log\mathcal{C}$ is empty (see above), and  it is known that for any linear normal
operator  in the Hilbert space the residual spectrum is empty.

(iv) By \cite[Theorem 1]{faa}
$$
\|R(\lambda,\log\mathcal{C})\|=\sup_{s\in \mathbb{R}}|\varphi(s)|=\frac{1}{d(\lambda)},
$$
where 
$$
d(\lambda)=\min_{s\in \mathbb{R}}|\lambda+\log(\frac12-\imath s)|=\mathrm{dist}(\lambda,\sigma(\log\mathcal{C}))
$$
(see the proof of (iii)). Further, it follows from (iii) that $\sigma(\log\mathcal{C})\cap \mathbb{R}=\{\log 2\}$. Therefore, each $\lambda\in \mathbb{R}, \lambda\ne \log 2$ is a regular point for 
$\log\mathcal{C}$, and  it is obvious  that $d(\lambda)=\lambda-\log 2$ for  $ \lambda> \log 2$. It remains to note that (iii) implies that $\lim_{\lambda\to-\infty}d(\lambda)=\pi/2$, since the lines $y=\pm\pi/2$ are horizontal asymptotes for the curve $\sigma(\log\mathcal{C})$ as $y\to\pm\pi/2$.
The proof is complete.
\end{proof}

Now we are in a position to study the inverse of the operator $\log\mathcal{C}$.

\begin{theorem}\label{th:inverse} The operator $\log\mathcal{C}$ has a bounded  inverse $(\log\mathcal{C})^{-1}$ that  belongs to the algebra $\mathcal{A}_u$. Moreover,  the scalar symbol of $(\log\mathcal{C})^{-1}$ is
  $$
\varphi(s)=  -\frac{1}{\log(\frac12-\imath s)}.
  $$

\end{theorem}

\begin{proof}  We put $L:=\log\mathcal{C}$ and $R_1:=R(1,\log\mathcal{C})$ for short. First note that $0$ and  $1$ are  regular points for $L$ by Theorem \ref{th:resolvent} (iii) above, and the assertion (ii) of this theorem
implies that $1$ is a regular point for $R_1$.  On the other hand, it is easy to verify that for any operator $L$
$$
L^{-1}=(I-L)^{-1}((I-L)^{-1}-I)^{-1}
$$
if all inverse operators exist. It follows that  
$$
(\log\mathcal{C})^{-1}=F(R_1),
$$
where $F(z)=z/(z-1)$ is a holomorphic function on the neighborhood of the spectrum of $R_1$ by Theorem \ref{th:resolvent} (ii), and $F(0)=0$. 
Since $R_1\in \mathcal{A}_u$  by Theorem \ref{th:resolvent} (i), we conclude that $(\log\mathcal{C})^{-1}$ belongs to $\mathcal{A}_u$ by 
\cite[Theorem 3.1]{LiflyandMirotin1}. Theorem \ref{th:resolvent} (i ) shows that operator $R_1$ satisfies all the conditions of \cite[Theorem 1]{faa}. By this theorem, this operator is unitary equivalent  to the
 operator of coordinate-wise multiplication by  its scalar symbol (see formula \eqref{varphi})
 $$
 \varphi_1(s)=\frac{1}{1+\log(\frac12-\imath s)}
 $$
  in the space $L^2(\mathbb{R}, \mathbb{C}^2)$  of $\mathbb{C}^2$-valued functions. Thus, the matrix symbol or $R_1$ is $\Phi_1={\rm diag}(\varphi_1,\varphi_1)$.
  It follows that  the matrix symbol or $F(R_1)=(\log\mathcal{C})^{-1}$ is $\Phi=F(\Phi_1)={\rm diag}(F(\varphi_1),F(\varphi_1))$ (see the proof of Theorem 3.1 in \cite{LiflyandMirotin1}).
  So, the scalar symbol of $(\log\mathcal{C})^{-1}$ is
  $$
\varphi(s)=  F(\varphi_1)(s)=-\frac{1}{\log(\frac12-\imath s)}.
  $$
\end{proof}

\section{Functions of the Ces\'aro operator in $L^2(\mathbb{R}_+)$}

In this section we briefly discuss how the results from the previous sections can be adopted to the  Ces\'aro operator in $L^2(\mathbb{R}_+)$.

\subsection{Holomorphic functions of the Ces\'aro operator in $L^2(\mathbb{R}_+)$}

We consider $L^2(\mathbb{R}_+)$ as a subspace of $L^2(\mathbb{R})$, and  for an operator $(A,D(A))$ in $L^2(\mathbb{R})$ we denote by $A_+$  the restriction $A|(D(A)\cap L^2(\mathbb{R}_+))$.
Note that  $ L^2(\mathbb{R}_+)$ is a closed invariant subspace for an operator  $\mathcal{H}_K$  if $K$ is supported in  $\mathbb{R}_+$ and  $K(u)|u|^{-\frac12}\in L^1(\mathbb{R})$.

Since $\sigma(\mathcal{C}_+)=\sigma(\mathcal{C})=\mathbb{T}+1$ \cite{BHSh}, \cite{faa} the expression $F(\mathcal{C}_+)$ makes sense for each function $F$  holomorphic in a neighborhood  of the circle $\mathbb{T}+1$. A preliminary lemma is required. 

\begin{lemma}\label{lm:1} Let $F$ be a holomorphic function in a neighborhood $N$ of the circle $\mathbb{T}+1$. Then
 $$
 F(\mathcal{C}_+)=F(\mathcal{C})_+.
 $$
\end{lemma}

\begin{proof} 
 One can assume that  $F(0)=0$. Let $\lambda\notin \sigma(\mathcal{C}_+)=\sigma(\mathcal{C})$.  For each $g\in L^2(\mathbb{R}_+)$ the equality
$f=R(\lambda,\mathcal{C})_+g$ means that $f\in  L^2(\mathbb{R})$ is the unique solution to the equation $(\lambda I-\mathcal{C})f=g$. Similarly,  the equality
$h=R(\lambda,\mathcal{C}_+)g$ means that $h\in  L^2(\mathbb{R}_+)\subset L^2(\mathbb{R})$ is the unique solution to the equation $(\lambda I-\mathcal{C})h=g$. Thus, $f=h$ and therefore $R(\lambda,\mathcal{C})_+=R(\lambda,\mathcal{C}_+)$. Now let $\Gamma$ denotes the boundary of any open annulus $U$ that contains  $\mathbb{T}+1$ such that $N$ contains its closure. Then
\begin{align*}
F(\mathcal{C})_+&:=\frac{1}{2\pi i}\left(\int_{\Gamma} F(\lambda)R(\lambda,\mathcal{C})d\lambda\right)_+\\
&=\frac{1}{2\pi i}\int_{\Gamma} F(\lambda)R(\lambda,\mathcal{C})_+d\lambda\\
&=\frac{1}{2\pi i}\int_{\Gamma} F(\lambda)R(\lambda,\mathcal{C}_+)d\lambda\\
&=F(\mathcal{C}_+).
\end{align*}
\end{proof}

As a simple corollary of Theorem \ref{th:1} and Lemma \ref{lm:1} we get the following result.

\begin{theorem}\label{th:3} Let $F$ be a holomorphic function in a neighborhood  of the circle $\mathbb{T}+1$, $F(0)=0$. Then  there is such function $K$, that  $F(\mathcal{C}_+)=\mathcal{H}_{K+}$, and conditions (a), (b), and (c) hold.

Conversely, if a function $K$   meets the conditions  (a), (b), and (c), where $F$ is  holomorphic  in some neighborhood  of  $\mathbb{T}+1$, then $F(0)=0$ and $\mathcal{H}_{K+}=F(\mathcal{C}_+)$. 
\end{theorem}

\begin{proof} Indeed, if $F$ is a holomorphic function in a neighborhood  of the circle $\mathbb{T}+1$, $F(0)=0$, then $F(\mathcal{C})=\mathcal{H}_{K}$ and   (a), (b), and (c) hold by Theorem \ref{th:1}. It follows by Lemma \ref{lm:1} that $F(\mathcal{C}_+)=F(\mathcal{C})_+=\mathcal{H}_{K+}$.

Conversely, if a function $K$   meets the conditions  (a), (b), and (c), where $F$ is  holomorphic  in some neighborhood  of  $\mathbb{T}+1$, then $F(0)=0$ and $F(\mathcal{C})=\mathcal{H}_{K}$ and therefore  $\mathcal{H}_{K+}=F(\mathcal{C})_+=F(\mathcal{C}_+)$.
\end{proof}

\begin{corollary}
Under the conditions of Theorem \ref{th:3} one has for the spectrum  $\sigma(\mathcal{H}_{K+})=F(\mathbb{T}+1)$.
\end{corollary}



\begin{example}\label{ex:P}
Let ${\rm Re}(\alpha)<\frac12$. Consider the  averaging operator of Boyd in $ L^2(\mathbb{R}_+)$ \cite{Boyd1}, \cite{Boyd}
$$
P_\alpha f(x)=x^{\alpha-1}\int_0^xt^{-\alpha}f(t)dt=\int_0^1u^{-\alpha}f(ux)\,du,
$$
where $P_0=\mathcal{C}_+$ is the continuous Ces\'{a}ro operator.    Then $P_\alpha=\mathcal{P}_{\alpha +}$ (see Example \ref{ex:1}). Therefore by Lemma \ref{lm:1}

  $$
P_\alpha=\alpha^{-1}((I-\alpha\mathcal{C}_+)^{-1}-I)
  $$

   It follows also, that the spectrum  $\sigma(P_\alpha)$ equals to the circle $F_\alpha(\mathbb{T}+1)$, where $F_\alpha$ is given by \eqref{Falpha}.
\end{example}

\begin{example}\label{ex:q}
Let ${\rm Re}(\alpha)<\frac12$. Consider the another   averaging operator of Boyd in $ L^2(\mathbb{R}_+)$ \cite{Boyd1}, \cite{Boyd}
$$
Q_\alpha f(x)=x^{-\alpha}\int_x^\infty t^{\alpha-1}f(t)dt=\int_1^\infty u^{\alpha-1}f(ux)\,du.
$$
   Then $Q_\alpha =\mathcal{Q}_{\alpha +}$ (see Example \ref{ex:2}), and by Lemma \ref{lm:1}

  $$
Q_\alpha=\mathcal{C}_+((I-\alpha)\mathcal{C}_+-I)^{-1}.
  $$

  We conclude also, that the spectrum  $\sigma(Q_\alpha)$ equals to the circle $F(\mathbb{T}+1)$, where $F$ is given by \eqref{F}.
\end{example}


\begin{example}\label{ex:5}
Let $\alpha>0$. Consider the generalized  Ces\'aro operator  of order $\alpha>0$ in $L^2(\mathbb{R}_+)$
$$
(\mathcal{C}_{\alpha +} f)(x)=\alpha x^{-\alpha}\int_0^x(x-t)^{\alpha-1}f(t)dt.
$$
Then $\mathcal{C}_{\alpha +}=\mathcal{H}_{K+}$ where $K(u)=\alpha(1-u)^{\alpha-1}\chi_{(0,1)}(u)$. In this case,
$(\mathcal{M}K)(z)=\alpha{\rm B}(\alpha,z)=F(z^{-1})$ where 
$$
F(z)=\frac{\Gamma(\alpha+1)\Gamma(z^{-1})}{\Gamma(\alpha+z^{-1})}.
$$
This function is  holomorphic in a neighborhood of $\mathbb{T}+1$ if and only if $\alpha\in\mathbb{N}$. So the operator $\mathcal{C}_{\alpha +}$     is  a  holomorphic function of  $\mathcal{C}_+$  if and only if $\alpha\in\mathbb{N}$. In this case $\mathcal{C}_{\alpha +}=F(\mathcal{C})$.
\end{example}

\subsection{Fractional powers of the Ces\'aro operator in $L^2(\mathbb{R}_+)$}

The operator   $\mathcal{C}_+$ is normal and the analog of Theorem \ref{th:2} is valid for $\mathcal{C}_+$, too. The case  $\alpha\in \mathbb{N}$ in \eqref{Salpha} below was proved by Boyd \cite{Boyd}. For another proof of Boyd's formula see \cite{LiflyandMirotin2}.

\begin{theorem}\label{th:4} Let $\mathrm{Re}(\alpha)>0$. Then
\begin{align}\label{Salpha}
(\mathcal{C}_+)^\alpha=H_{\alpha+}.
\end{align}
\end{theorem}

\begin{proof} Note that the operator  $\mathcal{C}$ has the form \eqref{HKa}      with  $a(u)=u\chi_{[0,1]}(u)$ and thus $a(u)\ge 0$ for all $u$. In this case the Corollary 4 in \cite{faa} is applicable. According to this corollary the operator  $\mathcal{C}_+$  is unitary equivalent to the operator $M_{\varphi_+^0}$ of multiplication by $\varphi_+^0$ in $L^2(\mathbb{R}_+)$. Thus, one can in fact to repeat the proof of Theorem \ref{th:2} simply replacing the operator $M_{\Phi_\mathcal{C}}$ by the operator $M_{\varphi_+^0}$.
\end{proof}

\begin{corollary}
Let $\log(\mathcal{C}_+)$ stands for the generator of the semigroup $(\mathcal{C}_+)^\alpha=H_{\alpha+}$. Then  $\log(\mathcal{C}_+)= (\log\mathcal{C})_+$.
\end{corollary}
This follows, e.g., from \cite[Subsection II.2.3, Corollary on p. 61]{EN}.

\begin{corollary}
Analogs of corollaries \ref{cor:5}---\ref{cor:8} and Theorem \ref{th:resolvent}  are valid for $H_{\alpha+}$ as well as for $H_{\alpha}$.
\end{corollary}

\section{Holomorphic functions of the Ces\'aro operator in $L^2[0,1]$}

We consider $L^2[0,1]$ as a subspace of $L^2(\mathbb{R})$, and  for an operator $A$ on $L^2(\mathbb{R})$ we denote by $A_1$  the bi-restriction (restriction and co-restriction) of $A$ to the space$L^2[0,1]$. 
This means that the domain of  $A_1$ is  $L^2[0,1]$ and  $A_1f=(Af)\chi_{(0,1)}$ for all  $f\in L^2[0,1]$.

Since $\sigma(\mathcal{C}_1)=\overline{\mathbb{D}}+1$ \cite{BHSh}, the expression $F(\mathcal{C}_1)$ makes sense for each function $F$ which is  holomorphic in a neighborhood  of the closed  disc $\overline{\mathbb{D}}+1$. Next, the following analog of Lemma \ref{lm:1} holds true.

\begin{lemma}\label{lm:2} Let $F$ be a holomorphic function in a neighborhood $N$ of the disc $\overline{\mathbb{D}}+1$. Then
 $$
 F(\mathcal{C}_1)=F(\mathcal{C})_1.
 $$
\end{lemma} 

The proof of this lemma is similar to the proof of Lemma \ref{lm:1}.

Now we have the following analog of Theorem \ref{th:3}.

\begin{theorem}\label{th:4} Let $F$ be a holomorphic function in a neighborhood  of the disk    $\overline{\mathbb{D}}+1$, $F(0)=0$. Then then there is such function $K$, that  $F(\mathcal{C}_1)=(\mathcal{H}_{K})_1$, and conditions (a), (b), and (c) hold.

Conversely, if a function $K$   meets the conditions  (a), (b), and (c), where $F$ is  holomorphic  in some neighborhood  of  $\overline{\mathbb{D}}+1$,  then $F(0)=0$ and $(\mathcal{H}_{K})_1=F(\mathcal{C}_1)$. 
\end{theorem}

\begin{proof} Indeed, if $F$ is a holomorphic function in a neighborhood  of the disc $\overline{\mathbb{D}}+1$, $F(0)=0$, then $F(\mathcal{C})=\mathcal{H}_{K}$ and   (a), (b), and (c) hold by Theorem \ref{th:1}. It follows by Lemma \ref{lm:2} that $F(\mathcal{C}_1)=F(\mathcal{C})_1=(\mathcal{H}_{K})_1$.

Conversely, if a function $K$   meets the conditions  (a), (b), and (c), where $F$ is  holomorphic  in some neighborhood  of  $\overline{\mathbb{D}}+1$, then  $F(0)=0$ as in the proof of Theorem \ref{th:1}, $F(\mathcal{C})=\mathcal{H}_{K}$ and therefore  $(\mathcal{H}_{K})_1=F(\mathcal{C})_1=F(\mathcal{C}_1)$.
\end{proof}

\begin{corollary} We have for the resolvent of $\mathcal{C}_1$
$$
R(\lambda,\mathcal{C}_1)=\lambda^{-2}(P_{\lambda^{-1}})_1+\lambda^{-1} I_1
$$
for all regular $\lambda$, i.e., for all  $\lambda$  exterior to the disc $\overline{\mathbb{D}}+1$.
\end{corollary}

\begin{proof} This follows from \eqref{R} and  Lemma \ref{lm:2}, since the spectrum of $\mathcal{C}_1$ equals  $\overline{\mathbb{D}}+1$.

\end{proof}


\end{document}